\documentclass[a4paper, 11pt, twoside, final]{article}

\newcommand{\fpageno}{3} %article first page number
\usepackage{graphicx}
\usepackage{amssymb,amsfonts,amsthm, amscd}
\usepackage{amsmath}
\usepackage{hyperref}
\usepackage[mathcal]{eucal}

\title{\normalsize ON A PROPERTY OF NORMING CONSTANTS OF\\
STURM--LIOUVILLE PROBLEM}

\author{\small {Yu.A. ASHRAFYAN and T.N. HARUTYUNYAN} \\ \small{Yerevan State University, Yerevan, Armenia}}

\newtheorem*{thm*}{Theorem A}
\newtheorem*{thm**}{Theorem B}

\def\a{\alpha}
\def\b{\beta}

\def\be{\begin{equation}}
\def\ee{\end{equation}}

\def\a{\alpha}
\def\v{\varphi}
\def\b{\beta}

\def\la{\lambda}

\def\l{\left(}
\def\r{\right)}
\def\lk{\left[}
\def\rk{\right]}

\def\il{\int\limits^{\pi}_0}
\def\mbc{{\mathbb C}}
\def\mbr{{\mathbb R}}

\begin{document}
\setcounter{page}{\fpageno} % first page number

\maketitle

\begin{abstract}
Using the Gelfand--Levitan method for the solution of inverse Sturm--Liouville problem we find a connection, which shows the dependence of norming constants on boundary conditions.
\vskip 7pt

\end{abstract}

\parindent=1cm

\section{Introduction.}\ Let $L(q, \a, \b)$ denote the Sturm--Liouville boundary value problem
\begin{gather}
\ell y\equiv -y''+q(x)y=\mu y,\quad x\in (0, \pi),\; \mu\in \mbc,\label{eq1}\\
y(0)\cos\a+y'(0)\sin\a=0,\quad \a\in (0, \pi),\label{eq2}\\
y(\pi)\cos\b +y'(\pi)\sin\b=0,\quad \b\in(0, \pi),\label{eq3}
\end{gather}
where $q$ is a real-valued, summable on $[0, \pi]$ function (we write $q\in L^1_{\mbr}[0, \pi]$).
By $L(q, \a, \ \b)$ we \ also \ denote \ the \ self-adjoint \ operator, \ generated \ by \ the \ problem \eqref{eq1}--\eqref{eq3} (see \cite{1}). It is known that under these conditions the spectra of the operator $L(q, \a, \b)$ is discrete and consists of real, simple eigenvalues \cite{1}, which we denote by $\mu_n=\mu_n(q,\a,\b)=\la_n^2(q, \a, \b)$, $n=0, 1, 2, \ldots$, emphasizing the dependence of $\mu_n$ on $q$, $\a$ and $\b$.

Let $\v(x,\mu,\a,q)$ and $\psi(x,\mu,\b,q)$ are the solutions of Eq. (1), which satisfy the initial conditions
$$\v(0,\mu,\a,q)=\sin\a,\   \ \v'(0,\mu,\a,q)=-\cos\a,$$

$$\psi(\pi,\mu,\b,q)=\sin\b, \   \ \psi'(\pi,\mu,\b,q)=-\cos\b,$$
correspondingly.
The eigenvalues $\mu_n=\mu_n(q, \a, \b)$, $n=0,1, 2, \ldots$, of $L(q, \a, \b)$ are the solutions of the equation
$$\Phi(\mu)=\Phi(\mu,\a,\b)\stackrel{def}{=}\v(\pi, \mu, \a)\cos \b+\v'(\pi, \mu,\a)\sin\b=0,$$
or the equation
$$\Psi(\mu)=\Psi(\mu,\a,\b)\stackrel{def}{=}\psi(0, \mu, \b)\cos\a+\psi'(0, \mu, \b)\sin \a=0.$$
According to the well-known Liouville formula, the wronskian $$W(x)=W(x,\v,\psi)= \v\cdot\psi'-\v'\psi$$ of the solutions $\v$ and $\psi$ is constant. It follows that $W(0)=W(\pi)$ and consequently $\Psi(\mu, \a, \b)=-\Phi(\mu, \a, \b).$
It is easy to see that the functions $\v_n(x)=\v(x, \mu_n, \a)$ and $\psi_n(x)=\psi(x,\mu_n, \b)$, $n=0, 1, 2, \ldots$, are the eigenfunctions, corresponding to the eigenvalue $\mu_n$.
Since all eigenvalues are simple, there exist constants $c_n=c_n(q, \a, \b)$, $n=0, 1, 2, \ldots$, such that
\begin{equation}\label{eq4}
\v(x,\mu_n)=c_n\cdot\psi(x,\mu_n).
\end{equation}
The squares of the $L^2$-norm of these eigenfunctions:
$$a_n=a_n(q,\a,\b)=\il |\v_n(x)|^2 dx,\    n=0, 1, 2, \ldots,$$

$$b_n=b_n(q,\a,\b)=\il |\psi_n(x)|^2 dx,\    n=0, 1, 2, \ldots$$
are called the norming constants.

In this paper we consider the case $\a,\b \in (0,\pi)$, i.e. we assume that $\sin \a \neq 0$ and
$\sin\b \neq 0$.
In this case we consider the solution $\tilde{\v}(x,\mu,\a,q):=\dfrac{\v(x,\mu,\a,q)}{\sin\a}$ of the equation (1) which has the initial values
$$\tilde{\v}(0,\mu,\a,q)=1, \ \tilde{\v}(x,\mu,\a,q)=-\cot\a,$$
and also we consider the solution $\tilde{\psi}(x,\mu,\b,q):=\dfrac{\psi(x,\mu,\b,q)}{\sin\b}$.
Of course, the functions $\tilde{\v}_n(x):=\tilde{\v}(x,\mu_n,\a,q)$ and $\tilde{\psi}_n(x):=\tilde{\psi}(x,\mu_n,\a,q)$, $n=0,1,2,\ldots,$ are the eigenfunctions, corresponding to the eigenvalue $\mu_n$. It follows from (4) that for norming constants
$\tilde{a}_n:=\|\tilde{\v}_n\|^2=\dfrac{a_n}{\sin^2\a},$\
$\tilde{b}_n:=\|\tilde{\psi}_n\|^2=\dfrac{b_n}{\sin^2\b}$
the following connections

\begin{equation}\label{eq5}
\tilde{b}_n=\frac{b_n}{\sin^2\b}=\frac{a_n}{c^2_n \sin^2\b}=
\frac{\tilde{a}_n \sin^2 \a}{c^2_n \sin^2\b}
\end{equation}

\noindent hold.

\section{The Main Result.}

The aim of this paper is to prove the following assertion.

\textbf{\textit {Theorem.}}
For the norming constants $\tilde{a}_n$ and $\tilde{b}_n$ the following connections  hold:

\begin{gather}
\ \frac{1}{\tilde{a}_0}-\frac{1}{\pi}+\sum_{n=1}^\infty\l\frac{1}{\tilde{a}_n}-\frac{2}{\pi}\r=\cot\a, \label{eq6}\\
\frac{1}{\tilde{b}_0}-\frac{1}{\pi}+\sum_{n=1}^\infty\l\frac{1}{\tilde{b}_n}-\frac{2}{\pi}\r=-\cot\b. \label{eq7}
\end{gather}

For the solution $\tilde{\v}$ it is well known the representation (see \cite{2,3})
\begin{equation}\label{eq8}
\tilde{\v}(x,\la,\a,q)=\cos{\la x}+\int\limits^x_0 G(x,t)\cos{\la t}dt,
\end{equation}
where for the kernel $G(x,t)$ we have (in particular) (see \cite{3})
\begin{equation}\label{eq9}
G(x,x)=-\cot\a + \frac{1}{2} \int\limits^x_0 q(s)ds.
\end{equation}
Besides, it is known that $G(x,t)$ satisfies to the Gelfand--Levitan integral equation
\begin{equation}\label{eq10}
G(x,t)+F(x,t)+\int\limits^x_0 G(x,s)F(s,t)ds=0,\  \ 0\leq t \leq x,
\end{equation}
where the function $F(x,t)$ is defined by the formula (see \cite{3})
\begin{equation}\label{11}
F(x,t)=\sum_{n=0}^\infty \l \frac{\cos{\la_n x}\cos{\la_n t}}{\tilde{a}_n}-\frac{\cos{n x}\cos{n t}}{a_n^0}\r,
\end{equation}
where $a_0^0=\pi$ and $a_n^0={\pi}/2$ for $n=1,2,\ldots$
It  easily follows from (9)--(11) that
$$G(0,0)=-F(0,0)=-\sum_{n=0}^\infty \l \frac{1}{\tilde{a}_n}-\frac{1}{a_n^0}\r=$$
\begin{equation}\label{eq12}
=-\l \frac{1}{\tilde{a}_0}-\frac{1}{\pi}\r - \sum_{n=1}^\infty \l \frac{1}{\tilde{a}_n}-\frac{2}{\pi}\r = -\cot\a.
\end{equation}
Thus, (6) is proved.

Let us now consider the functions ($n=0,1,2,\ldots$)
\begin{equation}\label{eq13}
p(x,\mu_n)=\frac{\v(\pi-x,\mu_n,\a,q)}{\v(\pi, \mu_n,\a,q)}=
\frac{\v(\pi-x,\mu_n)}{\v(\pi, \mu_n)}.
\end{equation}
Since $\v(x,\mu,\a,q)$ satisfies  the Eq. (1) and

$$p'(x,\mu_n)=-\frac{\v'(\pi-x,\mu_n)}{\v(\pi, \mu_n)}, \ \ p''(x,\mu_n)=\frac{\v''(\pi-x,\mu_n)}{\v(\pi, \mu_n)},$$
we can see, that $p(x,\mu_n)$ satisfies the equation

$$-p''(x,\mu_n)+q(\pi-x)p(x,\mu_n)=\mu_n p(x,\mu_n)$$
and the initial conditions
\begin{equation}\label{eq14}
p(0,\mu_n)=1, \ \ p'(0,\mu_n)=-\frac{\v'(\pi,\mu_n)}{\v(\pi,\mu_n)}=-(-\cot\b)=\cot\b=-\cot(\pi-\b).
\end{equation}
We also  have
$$p(\pi,\mu_n)=\frac{\v(0,\mu_n)}{\v(\pi,\mu_n)}=\frac{\sin\a}{\v(\pi,\mu_n)}=
\frac{\sin(\pi-\a)}{\v(\pi,\mu_n)},$$

$$p'(\pi,\mu_n)=-\frac{\v'(0,\mu_n)}{\v(\pi,\mu_n)}=-\frac{-\cos\a}{\v(\pi,\mu_n)}=
\frac{-\cos(\pi-\a)}{\v(\pi,\mu_n)}.$$
From this it follows that $p_n(x):=p(x,\mu_n)$ satisfies the boundary condition
$$p_n(\pi)\cos(\pi-\a)+p'_n(\pi)\sin(\pi-\a)=0, \ n=0,1,2,\ldots$$
Let us denote $q^{*}(x):=q(\pi-x)$. Since $\mu_n(q^{*},\pi-\b,\pi-\a)=\mu_n(q,\a,\b)$
(it is easy to prove and is well known, see  for example \cite{4}),
it follows, that $p_n(x),\ n=0,1,2,\ldots,$ are the eigenfunctions of the  problem
$L(q^{*},\pi-\b,\pi-\a)$, which have the initial conditions (14), i.e.
$p_n(x)=\tilde{\v}(x,\mu_n,\pi-\b,q^{*}), \ n=0,1,2,\ldots$

Thus, as in (12), for the norming constants $\hat{a}_n=\|p(\cdot, \mu_n)\|^2$ we have
\begin{equation}\label{eq15}
\l \frac{1}{\hat{a}_0}-\frac{1}{\pi}\r + \sum_{n=1}^\infty \l \frac{1}{\hat{a}_n}-
\frac{2}{\pi}\r = \cot(\pi-\b)=-\cot\b.
\end{equation}
On the other hand, for the norming constants $\hat{a}_n$, according to (4), (5) and (13), we have
$$\hat{a}_n=\il p^2(x,\mu_n)dx= \il \frac{\v^2(\pi-x,\mu_n)}{\v^2(\pi,\mu_n)}dx=$$
$$=-\frac{1}{\v^2(\pi,\mu_n)} \int\limits^0_{\pi} \v^2(s,\mu_n)ds=
\frac{1}{\v^2(\pi,\mu_n)} \int\limits^{\pi}_0 \v^2(s,\mu_n)ds=$$
$$=\frac{a_n(q,\a,\b)}{\v^2(\pi,\mu_n)}=
\frac{\tilde{a}_n \sin^2\a}{c^2_n \sin^2\b}=\tilde{b}_n.$$
Therefore, we can rewrite (15) in the form

$$\l \frac{1}{\tilde{b}_0}-\frac{1}{\pi} \r -
\sum_{n=1}^\infty \l \frac{1}{\tilde{b}_n}-\frac{2}{\pi}\r =-\cot(\pi-\b)= \cot\b.$$
Thus (7) is true and Theorem  is proved.\hfill $\Box$

\section{Remark.}\ It is known from the inverse Sturm--Liouville problems, that the set of eigenvalues
$\Big\{\mu_n\Big\}_{n=0}^\infty$ and the  norming constants $\Big\{\tilde{a}_n\Big\}_{n=0}^\infty$ uniquely determine the problem $L(q, \a, \b)$. That means, in particular, that we can determine $\{\tilde{b}_n\}_{n=0}^\infty$ by these two sequences. Now we will derive the precise formulae for these connections.

It is known that the specification of the spectra $\Big\{ \mu_n(q,\a,\b) \Big\}_{n=0} ^\infty $ uniquely determines the characteristic function $\Phi(\mu)$ (see  \cite{4}, Lemma 1($iii$); \cite{5},  Lemma 2.2)  and also its derivative $\partial \Phi(\mu)/ \partial \mu=\dot{\Phi}(\mu)$ (\cite5 Lemma 2.3).

In particular, if $\a,\b \in (0,\pi)$ the following formulas hold:
\begin{equation}\label{eq16}
\dot{\Phi}(\mu_0)=-\pi \sin\a \sin\b \prod_{k=1}^\infty \frac{\mu_k - \mu_0}{k^2}
\end{equation}
and (if $n \neq 0$, i.e. $ n=1,2,\ldots $)
\begin{equation}\label{eq17}
\dot{\Phi}(\mu_n)=-\frac{\pi}{n^2} \lk \mu_0 - \mu_n \rk \sin\a \sin\b \prod_{k=1, k \neq n}^\infty \frac{\mu_k - \mu_n}{k^2}.
\end{equation}
On the other hand, it is easy to prove the relation (see \cite{5}, Eq. (2.16) in Lemma 2.2 and see \cite{4}, Lemma 1 ($iii$))
\begin{equation}\label{eq18}
a_n=-c_n \cdot \dot{\Phi}(\mu_n).
\end{equation}

Taking into account the connections \eqref{eq5} and \eqref{eq16}--\eqref{eq18} we can find formulae for
$\dfrac{1}{\tilde{b}_0}$ \ and \
$\dfrac{1}{\tilde{b}_n}, \ n=1,2,\ldots:$

$$\frac{1}{\tilde{b}_0}=\frac{\tilde{a}_0}{\pi^2 \cdot  \l \displaystyle\prod_{k=1}^\infty \frac{\mu_k - \mu_n}{k^2}\r^2},$$

$$\frac{1}{\tilde{b}_n}=\frac{\tilde{a}_n n^4}{\pi^2 \cdot [\mu_0-\mu_n]^2 \cdot \l \displaystyle\prod_{k=1, k \neq n}^\infty \frac{\mu_k - \mu_n}{k^2}\r^2}.$$

So, we can change the second assertion in Theorem  by the following equation
$$\frac{\tilde{a}_0}{\pi^2 \cdot  \l \displaystyle\prod_{k=1}^\infty \frac{\mu_k - \mu_n}{k^2}\r^2}-\frac{1}{\pi}+$$
$$+\sum_{n=1}^\infty\l
\frac{\tilde{a}_n n^4}{\pi^2 \cdot [\mu_0-\mu_n]^2 \cdot \l \displaystyle\prod_{k=1, k \neq n}^\infty \frac{\mu_k - \mu_n}{k^2}\r^2}-
\frac{2}{\pi}\r=-\cot\b.$$

\

\

This   research   is  supported  by the Open Society Foundations--Armenia, within the Education program, grant  No  18742

%%
%% Here goes the Bibliography
%%

{\def\section*#1{}
\vskip 5pt
\begin{center}
\uppercase{\footnotesize{REFERENCES}}
\end{center}
%\vskip 2pt


\begin{thebibliography}{5}
\setlength{\itemsep}{-2mm}

\small


\bibitem{1}  \textbf{Naimark M.A.}  Linear Differential Operators.   M.: Nauka, 1969 (in Russian).
\bibitem{2}  \textbf{Gel'fand I.M., Levitan B.M.} On the Determination of a Differential Equation from its Spectral Function. Izv. Akad. Nauk. SSSR. Ser Mat., 1951, v. 15, p. 253--304  (in Russian).
\bibitem{3} \textbf{Yurko V.A.} An Introductions to the Theory of Inverse Spectral Problems. M.: Fizmatlit, 2007 (in Russian).
\bibitem{4}  \textbf{Isaacson E.L., Trubowitz E.} The Inverse Sturm--Liouville Problem, I. Com. Pure and Appl. Math., 1983, v. 36, p. 767--783
\bibitem{5} \textbf{Harutyunyan T.N.} Representation of the Norming Constants by Two Spectra. Electronic Journal of Differential Equations, 2010, 159, p. 1--10.
%etc.

\end{thebibliography}
\end{document}